\def\timestamp{%
Time-stamp: <brouwer_and_cardinals.tex: Thursday 06-07-2017 at 16:56:46 (cest)>}
\def\stripname Time-stamp: <#1 #2>{#2}
\edef\filedate{\expandafter\stripname\timestamp}

\documentclass[a4paper]{amsart}

\newenvironment{summ}[1]%
               {\par\medbreak\noindent\textbf{#1}. \ignorespaces}%
               {\par}

\newenvironment{comment}%
               {\par\smallskip\noindent\textsl{Comment}. \ignorespaces}%
               {\par}

\usepackage{amsrefs}
\newcommand\Zbl[1]{ Zentralblatt~#1}
\newcommand\ccard[1]{\overline{\overline{#1}}}
\DeclareMathSymbol\N0{AMSb}{`N}
\DeclareMathSymbol\Q0{AMSb}{`Q}
\begin{document}
\title{Brouwer and cardinalities}

\author[K. P. Hart]{Klaas Pieter Hart}
\address{Faculty of Electrical Engineering, Mathematics and Computer Science\\
         TU Delft\\
         Postbus 5031\\
         2600~GA {} Delft\\
         the Netherlands}
\email{k.p.hart@tudelft.nl}
\urladdr{http://fa.its.tudelft.nl/\~{}hart}



\date{\filedate}

\begin{abstract}
This note is a somewhat personal account of a paper that L.~E.~J.
Brouwer published
in~1908 and that dealt with the possible cardinalities of subsets of 
the continuum.
That paper is of interest because it represents the first time that Brouwer
presented his ideas on foundations in an international forum.

I found Brouwer's notions and arguments at times hard to grasp if
not occasionally perplexing.
I hope that this note contributes to a further discussion of the definitions
and reasonings as presented in Brouwer's paper.
\end{abstract}

\maketitle

\section*{Introduction}

A long time ago I needed a paper from the Proceedings of the IVth International
Congress of Mathematicians (1908)~\cite{AttiIV-II-cite} 
and I came across a paper by Brouwer with the title 
``Die Moeglichen Maechtigkeiten''~\cite{Brouwer1908A}, 
in those same proceedings.
Because of my interest in Set Theory I skimmed through it and came away with
the impression that its conclusion was that the set of countable ordinals does
not exist; this was news to me but I did not have time to study the paper
closely, so I made a copy and put it away for later.
`Later' turned out to be the occasion of this volume.

When I went back to the paper it turned out that my mind and memory had 
conflated a few sentences and that the last statement of the paper was not
something like 
``Von Cantors zweite Zahlenklasse kann deshalb kein Rede sein''
but
\begin{quote}
``Von anderen unendlichen M\"achtigkeiten als die abz\"ahlbare,
die abz\"ahlbar-unfertige und die continuierliche, kann gar keine Rede sein''.
\end{quote}
In fact, earlier in the paper Cantor's second number class, now known 
as the set of countably infinite ordinals, was declared to be not a set at all.
This piqued my interest because I have used~$\omega_1$ quite often in
my work; it looks and feels quite real to me.
I also thought, and still think, that the number of infinities is quite large,
certainly larger than~three.

So I started to read the paper more closely and I found it harder to understand
than expected; the definitions and arguments are not always as concrete as 
I would have liked them to be.
To some extent this can be said for many papers of that era 
(late nineteenth, early twentieth century): 
many notions had a meaning that was tacitly assumed to be the same to every 
reader.
At times though one finds definitions that are not more than synonyms or that
appeal to some sort of intuitive process that turns out to be very hard
to formalize.
An example of the latter can be found in Section~\ref{sec:what-going-on}:
Cantor's definition of the `cardinal number' of a set. 
As we shall see one cannot prove that this definition is formally sound; 
fortunately though, Cantor established a concrete and workable equivalent of
`having the same cardinal number' to prove his results.

One could say that in~1908, after finishing his dissertation,
Brouwer was where Cantor was in the~1870s: 
he was taking the first steps in an investigation of what one could 
achieve constructively and with the use of the actual infinite.
When one reads Brouwer's collected works, \cite{MR0532661}, it becomes 
apparent that it took some decades before things found their final form.
The continuum was later \emph{constructed\slash defined} using, basically, 
the partial order of intervals whose end points are rational and with a power 
of~$2$ as denominator.
This idea still lives on, for example in pointless topology and the theory
of continuous lattices.
It is remarkable that these steps ultimately lead to an impressive edifice
with many applications in constructive mathematics and computer science.

\subsection*{Outline}
Section~\ref{sec:reading} summarizes~\cite{Brouwer1908A}
and gives some comments that expand on the previous paragraphs.
I should emphasize that I read the paper first \emph{as it is}, without 
consulting any background material, because this appears to be the first time 
that Brouwer published some of his ideas on foundations in an international 
forum.
These two pages and ten lines would be what a reader at the time would have 
access to, unless they could read Dutch and could get their hand on 
Brouwer's thesis.
Of course my reading is influenced by my own mathematical experiences: 
I~fully embrace modern Set Theory, excluded middle, Axiom of Choice, and all.

In Section~\ref{sec:what-going-on} we turn to Brouwer's collected works,
and in particular his thesis, to get a better idea of where the material
in the paper is coming from.
To some extent I met the same problems as when reading the paper: 
it is not always clear what Brouwer actually meant; many definitions are 
not complete or even present.
There is, as promised above, a short aside on Cantor's definition of
cardinal number and the way he dealt with this definition.

The final section tries to see what the `standard' mathematical content 
of~\cite{Brouwer1908A} is.
It turns out that, 
Brouwer's assumptions of what `our mathematical intuition' is capable of aside, 
one can read the paper as a proof that closed subsets of the
real line are either countable or of the same cardinality as the real 
line itself. 

My own conclusion is that Brouwer's ideas were still under development. 
Indeed, when one reads the collected works, \cite{MR0532661}, and especially 
the notes, it becomes apparent that it took several years, if not decades,
before the notions that were presented in the paper in a vague way got their 
final definitions.

\section{Reading the paper}\label{sec:reading}

\begin{summ}{First paragraph}
The paper opens with a paragraph that attempts to describe or explain how
mathematical systems come about or are created.
It speaks of a primordial intuition%
\footnote{Wenn man untersucht, wie die mathematischen Systeme zustande kommen,
          findet man, dass sie aufgebaut sind aus der Ur-Intuition der 
          Zweieinigkeit.
          Die Intuitionen des continuierlichen und des discreten finden 
          sich hier zusammen, weil eben ein Zweites gedacht wird nicht 
          f\"ur sich, sondern unter Festhaltung der Erinnerung des Ersten.
          Das Erste und das Zweite werden also \emph{zusammengehalten}, und in
          dieser Zusammenhaltung besteht die Intuition des continuierlichen
          (continere = zusammenhalten).
          Diese mathematische Ur-Intuition ist nichts anders als die 
          inhaltslose Abstraction der Zeitempfindung, d.\,h.\ der Empfindung
          von `fest' und `schwindend' zusammen, oder von `bleibend'
          und `wechselnd' zusammen.}
of the `Zwei\-ein\-ig\-keit' (two-one-ness(?)).

Here the intuition of the continuous and the discrete come together, because
a Second is conceived for itself but while retaining the memory of 
the First.
The First and the Second are thus \emph{kept together}, and this keeping
together constitutes the intuition of the continuous.
This mathematical primordial~intuition is nothing but a meaningless
(inhaltlos)
abstraction of the experience of time \dots
\end{summ}

\begin{comment}
This was hard to make clear sense of; 
`primordial intuition', 
`Zwei\-ein\-ig\-keit',
`the First', 
`the Second', 
`the continuous' \dots, 
seem to lack concrete meaning and this makes it hard to appreciate
what this paragraph is intended to convey.
As mentioned above one can argue that the paper is a product of its time and 
that in this context these terms would make sense to a reader back in the day;
that may be true for some readers, but not for all, as the
review quoted in Section~\ref{sec:what-going-on} attests.
\end{comment}

\begin{summ}{Second paragraph: of the primordial intuition}
The first concrete information about the primordial intuition is given in the
next paragraph: it contains the possibility of the following developments.
\begin{enumerate}
\item The construction of the order type~$\omega$; when one thinks of
      the full primordial intuition as a new First, then one can think of
      a new Second that one calls `three', and so on.
\item The construction of the order type~$\eta$; when one considers
      the primordial intuition as a transition between `First itself' and
      the `Second itself' then the 'Interposition' has been created.
\end{enumerate}
(The symbols~$\omega$ and~$\eta$ are Cantor's notation for the order types 
 of the sets of natural and rational numbers respectively, 
 see~\cite{Cantor1895}.)
\end{summ}

\begin{comment}
Here we see two building steps: `take the next' and `insert between' and
it appears that Brouwer allows for infinitely many applications of these
and that he is willing to consider those infinitely many steps
as finished, so that the sets of natural and rational numbers 
(or rather sets that look like these)
are available.
Thus it appears that Brouwer accepts the actual infinite, although he
does not give a justification for this assertion about the primordial 
intuition; one who does not believe in the actual infinite may retort: 
``that is your intuition, not mine --- 
just saying `and so on' is not enough to create an infinite entity''.
Nowadays we would simply say ``we assume the Axiom of Infinity''
and then prove from this that $\omega$ and~$\eta$ exist.
\end{comment}

\begin{summ}{What more is possible}
The next two paragraphs start the discussion of possible powers
(M\"achtigkeiten).

Every mathematical system constructed using the primordial intuition can
itself be taken as a new unit and this explains the richness of the infinite
fullness of the mathematically possible systems, that however can all be
traced back to the two aforementioned order types.

When one looks at things in this way there would be only one infinite power,
the countably infinite and, indeed, other discrete systems than the countable 
finished ones cannot be built.
There are two ways in which it does make sense to consider higher powers
in mathematics.
\end{summ}

\begin{comment}
To me it is not quite clear what 
``taking a system as a new unit'' accomplishes.
One set-theoretic interpretation is that given a system, $S$~say, one can 
form~$\{S\}$; but this seems of limited use.
That this does not lead us out of the countable realm is clear: the closure
of the system that consists of $\omega$ and~$\eta$ under the 
map $x\mapsto\{x\}$ is countable.

The last sentence indicates that there is a place for the uncountable 
in mathematics after all.
\end{comment}

\begin{summ}{Countable unfinished}
The first kind of (mathematically) possible uncountable entities is described 
as follows. 

One can describe a method for building a mathematical system that creates 
from every given countable set that belongs to the system a new element 
of that system.
With such a method one can, as everywhere in Mathematics, only construct 
countable sets; the full system can never be built in this way because it 
cannot be countable.
It~is incorrect to call the whole system a mathematical set, for it is
not possible to build it finished, from the primordial intuition.

Examples: the whole of the numbers of the second number class, the whole of 
the definable points on the continuum, the whole of the mathematical systems. 
\end{summ}

\begin{comment}
Here we get a bit more information; up till now it was not quite clear, 
to~me anyway, how the building of structures was to proceed.
With hindsight one can interpret the `thinking of a Second given a First'
as `performing some kind of construction', though no general description of
what constitutes a construction is given.

From a modern point of view one can take issue with the statement that
`the full system' can never be built in this way. 
If we assume that the construction should have some concrete description,
say by means of formulas, then one can invoke the L\"owenheim-Skolem 
theorem and find a \emph{countable} structure with the property that things 
definable\slash constructible from its members belongs to the structure again.

The example of the second number class shows that, implicitly, Brouwer is 
allowing a bit more: when one builds the system of countable ordinals,
as Cantor did in~\cite{Cantor1883b}*{\S\,11}, then one needs two operations:
the first adds~$1$ to every ordinal to construct its successor;
if the set of ordinals constructed thus far has no maximum then the second
operation adds a smallest upper bound as a new ordinal.
Apparently Brouwer assumes that `the whole of ordinals constructed thus far'
belongs to the system.
In a countable elementary substructure of `the universe' the set of countable
ordinals that belong to that structure is itself not a member of that 
structure.
Thus Brouwer's assumption is more general; his construction steps may 
apparently also be applied to subsets that are not elements. 
The question then becomes: to which of the uncountably many subsets
are we allowed to apply these steps? 
\end{comment}

\begin{summ}{The continuum}
The second kind of (mathematically) possible uncountable entities is related
to the continuum. 

One can consider the continuum as a matrix of points or units and assume
that two points can be considered distinct if and only if their positions
can be distinguished on a certain scale of order type~$\eta$.%
\footnote{Man kann das mit dem discreten gleichberechtigten Continuum als 
          Matrix von Punkten oder Einheiten betrachten, und annehmen, 
          dass zwei Punkte dann und nur dann als verschieden zu betrachten
          sind, wenn sie sich in ihrer Lage auf einer gewissen skala 
          von Ordnungstypus~$\eta$ unterscheiden lassen.}
One then observes that the thus defined continuum will never let itself
be exhausted as a matrix of points, and one has to add the possibility
of overlaying a scale of order type~$\eta$ with a continuum to 
the method for building mathematical systems.%
\footnote{Man bemerkt dann, dass das in dieser Weise definierte Continuum sich
          niemals als Matrix von Punkten ersch\"opfen l\"asst, und hat der 
          Methode zum Aufbau mathematischer Systeme hinzugef\"ugt die 
          m\"oglichkeit, \"uber eine Skala vom Ordnungstypus~$\eta$ ein 
          Continuum (im jetzt beschr\"ankten Sinne) hinzulegen.}
\end{summ}

\begin{comment}
I had problems with these sentences and I reproduce the original German
in the footnotes.
I dropped ``mit dem discreten gleichberechtigten'' from the translation
because it would make the sentence too contorted; the words indicate that
the continuum is an equally valid notion as the discrete.
There are three terms that are new: the continuum, matrix, and scale.
Most likely the continuum means the real line or an interval thereof.
As to matrix, one of the definitions in my dictionary, \cite{Chambers:Rev13}, 
is `the place in which anything is developed or formed'; we can only guess,
in this case, how this developing or forming is supposed to happen.
Finally, scale; since the continuum is linearly ordered this probably refers
to a subset whose order~type is equal to~$\eta$ and two points are distinct
if there is a point from the scale between them.

The second sentence states that the continuum-as-matrix-of-points 
can never be exhausted, but gives no reason, and tells us that this 
leads to a construction method: every scale of order type~$\eta$ 
can be covered by or completed to a continuum.

Also, in the first sentence the continuum appears to be a given object 
and in the second sentence it is suddenly `thus defined', without any 
recognizable intervening definition.
\end{comment}

\begin{summ}{Arbitrary subsets of the continuum}
Now begins an investigation into the possible cardinalities of subsets
of the continuum.

Let now $M$ be an arbitrary subset of some set of the power of the continuum.
Then one can map this continuum onto a linear continuum between~$0$ and~$1$
and thus this subset appears linearly ordered in this continuum.
In the set~$M$ there is a countable subset~$M_1$, 
with whose help $M$~is to be defined.
This set~$M_1$ can in the following manner be related to the scale of 
numbers~$a/2^n$ between~$0$ and~$1$.

We approximate the aggregate of the points of $M_1$ using binary fractions.
Every digit is either determined by its predecessor, or the choice
between~$0$ and~$1$ is still open;
we construct a branching system, where each branch continues in one direction
if the choice is not free, and splits itself when it is free.
After this we destroy each branch that does not split from some moment on 
from the first moment after which no split occurs.
When this is done to all such branches then we apply this procedure to the
remaining system.
Ultimately we will be left with an empty system or an infinite one in which 
every branch keeps splitting.
In the latter case $M_1$~has subsets of order type~$\eta$, in the former 
case it does not.

It may be that besides the countable set~$M_1$, of which it is stated that
it belongs to~$M$, the definition of~$M$ requires the specification
of a second countable set~$M_2$ that certainly must not belong to~$M$.
However, after these sets have been set up the determination of~$M$
can be completed in only one way, namely, in case $M_1$ has subsets of 
type~$\eta$ by executing the operation of continuous-making in one
or more of these sets while, of course, deleting the points excluded
by~$M_2$.

If this operation is executed at least once then the power of~$M$ is
that of the continuum, otherwise $M$~is countable.
\end{summ}

\begin{comment}
The first paragraph starts out a bit odd: 
$M$~is a subset of a set of the power of the continuum and in the same 
sentence that set is already a continuum itself 
(in the German original the first two sentences are one).
The next step is to transfer $M$ to the unit interval. 
The statement after that is rather vague: what does it mean that $M$~is 
to be defined with the help of~$M_1$, and why should it necessarily be a 
subset?

The description of the way in which the points of~$M_1$ are related to
binary fractions (numerator odd, denominator a power of~$2$) is a bit iffy
--- it speaks of individual digits while dealing with all of~$M_1$ ---
but basically sound; one can see the resulting branching system as the set 
of intervals of the form $[a2^{-n},{(a+1)2^{-n}}]$ that intersect~$M_1$,
ordered by inclusion (so as a tree it grows upside down).
The remainder of the second paragraph describes the standard 
Cantor-Bendixson procedure applied to the tree, rather than the set~$M_1$.

The beginning of the third paragraph is a bit mystifying in that it is not
specified how $M_2$~is instrumental in the making of~$M$ and there is no reason
why it should not belong to~$M$ 
(most likely the `not belong to' means `is disjoint from').
Finally, why should $M$ be the result of applying the method of continuous
making to~$M_1$ or its subsets of type~$\eta$, and the deletion of the points
of~$M_2$?
\end{comment}

\begin{summ}{The conclusion}
Thus there exists just one power for mathematical infinite sets, to wit the
countable.
One can add to these:
\begin{enumerate}
\item \emph{the countable-unfinished}, but by this is meant a \emph{method}, 
      not a set
\item \emph{the continuous}, by this is certainly meant something finished,
      but only as \emph{matrix}, not as a set.
\end{enumerate}
Of other infinite powers, besides the countable, the countable-unfinished, 
and the continuous, one cannot speak.
\end{summ}

\begin{comment}
Taken at face value this does not really explain anything; in both cases
the notion is defined in terms of an~other, undefined, notion.
What \emph{is} a method? 
What \emph{is} a matrix and how does it work?  

Another point is that above $M$~is assumed to be a subset of the continuum.
At the end it turns out that $M$ may be countable or that it contains
something that is a bijective image of a continuum, presumably this 
disqualifies it from being a set but there is no explicit mention of this.
\end{comment}

\section{What is going on?}
\label{sec:what-going-on}

To repeat, the previous section was devoted to the paper~\cite{Brouwer1908A},
nothing more, nothing less.
To me the writing makes it hard to appreciate the ideas that it wants 
to express.
It turns out I was not the only one who was mystified.
Here is a review of the paper from 
\textsl{Jahrbuch \"uber die Fortschritte der Mathematik}
(available on the Zentralblatt website).
\begin{quote}
Referent bekennt, da\ss{} die Betrachtungen des Verfassers ihm nicht 
v\"ollig klar sind; er beschr\"ankt sich also darauf, seine Schl\"usse 
w\"ort\-lich abzuschreiben: 
``Es existiert nur eine M\"ach\-tig\-keit f\"ur 
  mathematische unendliche Mengen, n\"amlich die abz\"ahlbare. 
  Man kann aber hinzuf\"ugen: 
  1.~die abz\"ahlbar-unfertige, aber dann wird eine Methode, keine 
     Menge gemeint;
  2.~die kontinuierliche, dann wird freilich etwas Fertiges gemeint, aber nur 
    als Matrix, nicht als Menge. 
  Von anderen unendlichen M\"achtigkeiten, ab der abz\"ahlbaren, der 
  ab\-z\"ahl\-bar-unfertigen, und der kontinuierlichen, kann gar keine 
  Rede sein.''
\end{quote}

To see what I may have missed I turned to Brouwer's collected 
works~\cite{MR0532661}, where \cite{Brouwer1908A}~is reproduced on 
pages~102--104.
In the notes we learn that the paper is a summary of a section
in chapter~1 of Brouwer's thesis~\cite{Brouwer1907}, which 
can be found in English translation in~\cite{MR0532661}.

That chapter opens with an inventory of what the mathematical intuition
may take for granted.
First there are the number systems that we all know:
natural numbers, integers, fractions, irrational numbers.
The construction of the latter is by means of Dedekind cuts, not all at once
but step-wise, although it is not quite clear to me how and when certain 
numbers can be considered known or constructed.
As in~\cite{Brouwer1908A} the totality of known numbers is at any point 
in the development still countable.
Next is the continuum, which is also taken as given and its description
is not very concrete but from the things that Brouwer does with it 
it~becomes apparent that one should, as above, have the real line in mind.
There is a construction of a dense copy of the dyadic rational numbers in the
continuum and, in a drastic change of pace, addition and multiplication
are defined and characterized using one- and two-parameter 
transformation groups acting on the continuum.
In~about 52~pages one then finds how to build various types of geometries
out of the continuum.
I~sketched this so that we know what went before when one comes to the section 
that~\cite{Brouwer1908A} was summarizing.

With the notation as in Section~\ref{sec:reading} this section of the thesis 
makes the relation between $M$, $M_1$, and~$M_2$ somewhat more explicit: 
$M_1$~should be dense in~$M$ so that every other point of~$M$ can be 
approximated by points of~$M_1$; the set~$M_2$ is subtracted from~$M$.
Rather than simply stating that $M$~has the same power as the continuum
when $M_1$~has subsets of order type~$\eta$ there is an indication
of how $M$~can be mapped onto the continuum.

What did not make it into~\cite{Brouwer1908A} is the conclusion that
``it appears that this solves the continuum problem, by adhering strictly 
to the insight: one cannot speak of the continuum as a point~set other than
in relation to a scale of order type~$\eta$''.

\subsection*{A look through modern eyes: Cantor}

It has been suggested to me that I~should take the time of writing into 
consideration and to look on the Brouwer of~\cite{Brouwer1907}
and~\cite{Brouwer1908A} as a nineteenth-century mathematician.
And it must be said that many papers from that era tacitly assume that 
many things are understood by everyone to mean the same thing, even 
though no formal definition seems available. 

However sometimes such an appeal to `common knowledge' or a `common intuition'
can lead one astray.
As a case in point consider the following definition by Georg Cantor, 
in~\cite{Cantor1895}, of cardinal number:
\begin{quote}
,,M\"achtigkeit`` oder ,,Cardinalzahl`` von $M$ nennen wir den Allgemeinbegriff,
welcher mit Hilfe unseres activen Denkverm\"ogens dadurch aus der Menge~$M$
hervorgeht, da\ss{} von der Beschaffenheit ihrer verschiedenen Elemente~$m$
und von der Ordnung ihres Gegebenseins abstrahirt wird.  

Das Resultat dieses zweifachen Abstractionsakts, die Kar\-di\-nal\-zahl oder
M\"achtig\-keit von~$M$, bezeichnen wir mit
$$
\ccard{M}
$$
\end{quote}
Just like Brouwer's first paragraph in~\cite{Brouwer1908A} this sounds
wonderful in German but it does not quite deliver on its promise.
It does not specify how this abstraction should be carried out.
In fact one can argue that it cannot be carried out.

When one reads the first few pages of~\cite{Cantor1895}, where Cantor develops
this notion of cardinality, one will see that Cantor thinks of some process
whereby a set~$M$ morphs into a `standard' set~$\ccard{M}$. 
This process has the property that there is a bijection between~$M$ and~$N$
if and only if $\ccard{M}=\ccard{N}$ and the if-direction is established
via bijections between $M$~and~$\ccard{M}$, and between $N$~and~$\ccard{N}$. 
However, Theorem~11.3 in~\cite{MR0396271}
(attributed to Pincus~\cite{Pincus1969}) states that it is consistent with 
Zermelo-Fraenkel Set~Theory that no such assignment $M\mapsto\ccard{M}$ exists.

By contrast, the Axiom of Choice does enable one to construct a class 
of `standard' sets (called cardinal numbers) against which all sets can be 
measured and, indeed, an assignment $M\mapsto\ccard{M}$ with the properties
desired by Cantor.
Those constructions involve some non-trivial work and the resulting assignment
is most certainly not the product of some kind of `Act of Abstraction'.

Fortunately, whenever Cantor proved that two sets had the same cardinal number
he did so by providing explicit bijections between them. 
But this came at the end of two decades of work on sets and by then it 
was quite clear how to make the intuition explicit and how to actually work 
with cardinal numbers: use injections and bijections to establish inequalities
and equalities respectively.

\subsection*{A look through modern eyes: Brouwer}

To return to \cite{Brouwer1908A}, I found it frustrating reading simply 
because, as I mentioned before, there is apart from 
``overlaying a scale of order type~$\eta$ with a continuum'' 
no description of what construction steps are available (or allowed).

We know from later developments that Brouwer intended to be as 
explicit\slash constructive as possible.
But, at least in~\cites{Brouwer1907,Brouwer1908A}, there is a step that,
with today's knowledge, is not constructive at all.

The choice, or construction, of a scale of order type~$\eta$ in the continuum
is typical of proofs in Analysis that involve sequences:
generally one proves that individual choices are possible and then uses
words like ``and so on'' to convince us that we can actually construct a 
\emph{sequence} of choices.
That of course requires an instance of (at least) the Countable
Axiom of Choice, and one can create models of Set~Theory 
with infinite sets that have no injective map from the set of natural numbers
into them.
 
In fact, one can create a model in which there is an ordered continuum
whose subsets \emph{that are in the model} are finite unions of intervals;
thus, in that model that ordered continuum does not have a scale of
order type~$\eta$.
The easiest construction takes the unit interval~$I$ and the group~$G$ of 
increasing bijections from~$I$ to itself with finitely many fixed points.
The permutation model with $I$ as its set of atoms and determined by the
filter of subgroups of the form $G_F=\{g\in G:(\forall x\in F)(gx=x)\}$
is as required, compare~\cite{MR0396271}*{Chapter~4}.

\section{Interpretation}
\label{sec:interpretation}

Let us try to make `standard' sense of Brouwer's argument for the
possible powers of subsets of the continuum.

When showing that $M$ is countable or has the power of the continuum
Brouwer invokes the step of completing a set of order type~$\eta$ 
to a continuum.
The paper is not very explicit as to how that is supposed to happen
when that set is nowhere dense in some other continuum.
The thesis gives a bit more detail: a description of how
``one continuum may cover another one with gaps''.

Consider for example the set~$E$ of points in the unit interval~$I$ with a 
finite binary expansion whose every \emph{odd-numbered} digit is equal to~$0$; 
this set has order type~$\eta$ and so it can be completed to a continuum.
We can take this continuum, call it~$K$, to be the unit interval itself, upon 
identifying~$E$ with~$D$, the set of \emph{all} numbers with a finite binary 
expansion.
Simply map a point $\sum_{i=1}^nd_i2^{-i}$ of~$D$ to the point 
$\sum_{i=1}^nd_i2^{-2i}$ of~$E$.
Now, $E$ is a subset of both~$K$ and~$I$, dense in the former and nowhere
dense in the latter; the problem then is to relate~$K$ to the subset 
that $E$ determines when considered as a subset of~$I$.

In the thesis Brouwer stipulates that, given a point~$x$ in the unit interval,
we transfer its approximation sequence from~$D$ to~$E$ and then
take the limit of that transferred sequence as the point associated to~$x$.
The problem, for me, lies with the `its approximating sequence'; each point
can be approximated by many sequences, from below and from above.
For points like the $\frac12$ in~$K$ this creates an ambiguity: should we 
associate it to~$\frac14$ or to~$\frac1{12}$ in~$I$?

\subsection*{What sets do actually exist?}

Rather than dwell on this ambiguity let us see if we can read between the lines
and describe the subsets of the unit interval that Brouwer would consider 
constructible.

The steps described in the thesis and in the paper enable us to construct
the following types of sets.
\begin{enumerate}
\item countable
\item closed
\item unions and differences of sets of the first two types
\end{enumerate}
The third step is clear: if one can make two sets then one can also make their
union and difference.
That the countable subsets of the unit interval are constructible is 
more or less true by definition: an enumeration of the set specifies it.
At every point in such a definition one can insert `constructive' 
or `definable' if one wants to delineate the type of enumeration being used.

To see that Brouwer's steps can also yield closed sets consider a countable
countable subset~$M_1$ of the unit interval that has been constructed
already.
We show how its closure~$M$ may be constructed.

The branches of the tree to which Brouwer applies his reduction procedure 
determine points of the closure of~$M_1$ and vice versa.
To every deleted branch~$b$ there corresponds a node~$s_b$ in the tree:
the place where $b$~is pruned from the tree.
This assignment $b\mapsto s_b$ is one-to-one, hence in total only countably
many branches will be pruned.

So if the pruning procedure leaves the tree empty then the closure of~$M_1$
is countable.

If it leaves a non-empty tree then the points of~$M_1$
that correspond to branches in this tree form a dense-in-itself subset~$J$
of the unit interval.
It need not have order type~$\eta$ because there may be points in~$J$ that are 
end~points of some maximal interval in the complement of the closure of~$M_1$.
However, these are the only exceptions, so that removing the set~$J_2$
of all left end~points from~$J$ will yield a set~$J_1$ of type~$\eta$.
We can complete~$J_1$ to a continuum~$K$.
Although the map from~$K$ to the unit interval suggested by Brouwer requires
a choice whenever a point of~$K$ represents a complementary interval
of~$M$ we can opt to choose the right end~point.
In this way the image of the map is the closure of~$J_1$ minus the set~$J_2$.
This constructs $M\setminus J_2$ and taking its union with~$J_2$
will yield the set~$M$.

\subsubsection*{Comments}

The arguments given above are a translation of those of Brouwer to our 
standard situation.
In later papers Brouwer became more explicit in his 
description\slash construction of the continuum; I recommend the notes
in~\cite{MR0532661} as a road map through the papers that deal with
intuitionistic set theory.
It appears that our translation remains valid under the new definitions.

One can also interpret Brouwer's arguments as an~other proof of Cantor's
theorem from~\cite{Cantor1884}*{\S\,19} that closed subsets of the unit 
interval either are countable or have the same cardinality as the unit 
interval itself.
If one does not worry about definability issues then every closed subset
of the continuum can occur as an~$M$: simply take $M_1$ to be a countable 
dense subset of~$M$, consider it constructible and obtain $M$ from it
as its closure.
In case $M$~is uncountable Brouwer's procedure yields an injective map
from the unit interval into~$M$.
This is different from Cantor's proof, which produced a surjection
from~$M$ onto the unit interval.

Also, the statement about Cantor's second number class was not a complete
figment of my imagination: one can find it in~\cite{Brouwer1907}
as number~XIII in a list of propositions 
``to be defended together with the thesis''.
It reads ``De tweede getalklasse van \textsc{Cantor} bestaat niet''
(The second number class of \textsc{Cantor} does not exist). 

\bigskip
I would like to thank, in chronological order,
Teun Koetsier, Jan van~Mill, and Jan~Willem Klop
for comments on earlier versions of this note; these led to noticeable 
improvements.

\smallskip
I end with a caveat:
on the 2\textsuperscript{nd} of September 2016,
at a meeting on the occasion of the transfer of the Brouwer archive to
the \textsl{Noord-Holland Archief},
J. Korevaar related an anecdote about Brouwer, see~\cite{Korevaar2016}.
At one of the monthly meetings of the Dutch Mathematical Society in the
1950s some mathematicians were explaining to each other what Brouwer
had meant in some paper.
Unnoticed by everyone Brouwer had entered the room and after a while
he ran to the board exclaiming: 
$$
\hbox{``You have all misunderstood!''}
$$

\end{document}